\documentclass[12pt]{article}
\usepackage{amsmath}
\usepackage{amssymb}
\font\tencmmib=cmmib10 \skewchar\tencmmib '60
\newfam\cmmibfam
\textfont\cmmibfam=\tencmmib

\def\bbox{\quad\hbox{\vrule \vbox{\hrule \vskip2pt \hbox{\hskip2pt
\vbox{\hsize=1pt}\hskip2pt} \vskip2pt\hrule}\vrule}}
\def\lessim{\ \lower4pt\hbox{$
\buildrel{\displaystyle <}\over\sim$}\ }
\def\gessim{\ \lower4pt\hbox{$\buildrel{\displaystyle >}
\over\sim$}\ }

\def\la{{\Bigl\langle}}
\def\ra{{\Bigr\rangle}}

\def\qed{\hfill\break\rightline{$\bbox$}}
\newcommand{\e}{\mathbb{E}}
\newcommand{\p}{\mathbb{P}}

\newcommand{\Reals}{\mathbb{R}}

\newtheorem{proposition}{Proposition}
\newtheorem{lemma}{Lemma}
\newtheorem{theorem}{Theorem}

\makeatletter
\@addtoreset{equation}{section}

\makeatother

\font\tencmmib=cmmib10 \skewchar\tencmmib '60
\newfam\cmmibfam
\textfont\cmmibfam=\tencmmib

\def\bbox{\quad\hbox{\vrule \vbox{\hrule \vskip2pt \hbox{\hskip2pt
\vbox{\hsize=1pt}\hskip2pt} \vskip2pt\hrule}\vrule}}
\def\lessim{\ \lower4pt\hbox{$
\buildrel{\displaystyle <}\over\sim$}\ }
\def\gessim{\ \lower4pt\hbox{$\buildrel{\displaystyle >}
\over\sim$}\ }

\def\go0{\to 0}

\def\la{\langle}

\def\leftitem#1{\item{\hbox to\parindent{\enspace#1\hfill}}}

\def\qed{\hfill\break\rightline{$\bbox$}}

\def\ra{\rangle}

\def\sg{\sigma}

\def\sg2{\sigma^2}

\def\__{_{\infty}}

\begin{document}

\title{
Ghirlanda-Guerra identities and ultrametricity: An elementary proof in the discrete case.}

\author{ 
Dmitry Panchenko\thanks{
Department of Mathematics, Texas A\&M University, 
Mailstop 3386, College Station, TX, 77843,
email: panchenk@math.tamu.edu. Partially supported by NSF grant.}\\
{\it Texas A\&M University}
}
\date{}

\maketitle

\begin{abstract}
In this paper we give another proof of the fact  that a random overlap array, which satisfies the Ghirlanda-Guerra identities 
and whose elements take values in a finite set, is ultrametric with probability one. The new proof bypasses random change 
of density invariance principles for directing measures of such arrays and, in addition to the Dobvysh-Sudakov representation, 
is based only on elementary algebraic consequences of the Ghirlanda-Guerra identities.
\end{abstract}
\vspace{0.5cm}



\section{Introduction and main result.}
In this paper we will give a  simplified proof of the main result in \cite{PGG}.
Let us consider an infinite random array $R= (R_{l,l'})_{l,l'\geq 1}$ which is symmetric, 
non-negative definite and weakly exchangeable, which means that for any $n\geq 1$ and for any
permutation $\rho$ of $\{1,\ldots, n\}$ the matrix $(R_{\rho(l),\rho(l')})_{ l,l' \leq n}$
has the same distribution as $(R_{l,l'})_{ l,l' \leq n}.$ 
We assume that diagonal elements $R_{l,l}=1$ and non-diagonal elements 
take finitely many values, 
\begin{equation}
\p\bigl(R_{1,2}=q_l\bigr) = p_l
\label{finite}
\end{equation}
for some  $-1 \leq q_1< q_2 <\ldots < q_k\leq 1$ and $p_l>0, p_1+\ldots+p_k = 1.$
The array $R$ is said to satisfy the Ghirlanda-Guerra identities \cite{GG} if for any $n\geq 2$ 
and any bounded measurable functions $f = f\bigl((R_{l,l'})_{1\leq l,l' \leq n}\bigr)$
and  $\psi:\Reals\to\Reals,$
\begin{equation}
\e f \psi(R_{1,n+1}) = \frac{1}{n}\, \e f\, \e \psi(R_{1,2}) + 
\frac{1}{n} \sum_{l=2}^{n} \e f \psi(R_{1,l}).
\label{GG}
\end{equation}
By the positivity principle of Talagrand  (\cite{SG}, \cite{SG2}), the Ghirlanda-Guerra identities imply that  
$R_{1,2}\geq 0$ with probability one and, therefore,  we can assume that $q_1\geq 0$.  

\begin{theorem}\label{ThMain}(\cite{PGG})
Under assumptions (\ref{finite}) and (\ref{GG}), the array $R$ is ultrametric, 
\begin{equation}
\p\bigl(R_{2,3}\geq \min(R_{1,2},R_{1,3})\bigr)=1.
\label{ultra}
\end{equation}
\end{theorem}
Another way to express the event in (\ref{ultra}) is to say that 
\begin{equation}
R_{1,2} \geq q_{l}, 
R_{1,3} \geq q_{l}
\Longrightarrow
R_{2,3} \geq q_{l}
\,\,\mbox{ for all }\,\, 
1\leq l\leq k.
\label{ultraeq}
\end{equation}
Infinite arrays that satisfy the Ghirlanda-Guerra identities arise as the limits of the overlap arrays  in the Sherrington-Kirkpatrick  
spin glass models (see e.g. \cite{SG2}, \cite{PGGmixed}). The assumption (\ref{finite}) is purely technical (and unfortunately is not satsified in the most important situations). The first ultrametricity result 
was proved in \cite{AA} under different conditions which also included (\ref{finite}), but instead of (\ref{GG}) the authors
worked with the Aizenman-Contucci stochastic stability \cite{AC}. The original proof of Theorem \ref{ThMain} in \cite{PGG} utilized
a key idea from \cite{AA}, namely,  the existence of directing measures guaranteed by the Dovbysh-Sudakov representation 
result in \cite{DS}, and we will still rely on this representation here. However, we will completely avoid proving any
invariance principles under random changes of density for the directing measure, which played crucial roles both in
\cite{AA} and \cite{PGG} and our new induction will be quite elementary in nature.  M. Talagrand gave a proof of  Theorem \ref{ThMain} in \cite{Tal-New}  that did not use the Dovbysh-Sudakov representation but  still used the invariance 
principle from \cite{PGG}.
The Dovbysh-Sudakov representation  \cite{DS}  (for detailed proof see \cite{PDS}) states that given
a symmetric, non-negative definite and weakly exchangeable array $R,$ there exists a random measure $\mu$ on 
$H\times[0,\infty),$ where $H$ is a separable Hilbert space, such that $R$ is equal in distribution to the array\begin{equation}
\bigl( \sigma^l\cdot \sigma^{l'} + a^l\,\delta_{l,l'}\bigr)_{l,l'\geq 1}
\label{repr}
\end{equation}
where $(\sigma^l, a^l)$ is an i.i.d. sequence from $\mu$ and $\sigma\cdot \sigma'$ denotes the scalar product  on $H$.
Let us denote by $G$ the marginal of $\mu$ on $H.$ The following simple consequence of the Ghirlanda-Guerra
identities (\ref{GG}) was proved in Theorem 2 in \cite{PGG}.
 \begin{proposition}\label{Th2}
 Under (\ref{finite}) and (\ref{GG}), the random measure $G$ is (countably) discrete and is concentrated on the sphere 
 of radius $\sqrt{q_k}$ with probability one.
 \end{proposition} 
In particular, this implies that $a^l=1-q_k$ in (\ref{repr}) and without loss of generality we can redefine the array by
$R_{l,l'} = \sigma^l\cdot\sigma^{l'}$
for an i.i.d. sequence $(\sigma^l)$ from $G.$ Since $R_{l,l'}=q_{k}$ if and only if $\sigma^l=\sigma^{l'}$,
we have
$$
\p(R_{1,2}=q_{k}, R_{1,3} = q_{k}, R_{2,3}<q_{k}) = 0,
$$
which proves ``ultrametricity at the level $k$" in the sense of (\ref{ultraeq}). As in \cite{AA} and \cite{PGG}, we would like to find a way
to make an induction step and prove ``ultrametricity at the level $k-1$". The main new idea of the paper will be to consider
the distribution of the array $(R_{l,l'})$ conditionally on the event that all replicas $(\sigma^l)$ are different and prove that this
new distribution is well-defined and satisfies all the conditions of the Dovbysh-Sudakov representation. Since on the above event 
the elements of the new array can not take value $q_k$, the induction step will follow.

\section{Proof.}

By Proposition \ref{Th2}, $G = \sum_{l\geq 1} w_l \delta_{\xi_l}$ for some random weights $(w_l)$ and random sequence 
$(\xi_l)$ in  $H$ such that $\xi_l \cdot \xi_l=q_k.$ Let us denote by $\la\cdot\ra$ the average with respect to $G^{\otimes \infty}$ 
and by $\e$ the expectation with respect to the randomness of $G.$ With these notations, the Ghirlanda-Guerra identities 
(\ref{GG}) can be rewritten as
\begin{equation}
\e \la f_n \psi(R_{1,n+1}) \ra  = \frac{1}{n} \e \la f_n\ra \, \e\la \psi(R_{1,2}) \ra + 
\frac{1}{n}\sum_{l=2}^{n} \e \la f_n \psi(R_{1,l})\ra.
\label{GG1}
\end{equation}
For each $n\geq 2$, let us consider the event
\begin{equation}
A_n = \{R_{l,l'} \not = q_k, \, \forall 1\leq l<l'\leq n\}
\end{equation}
and let $\p_n$ be the distribution of the $n\times n$ matrix  $R^n = (\sigma^l\cdot \sigma^{l'})_{l,l'\leq n}$ conditionally on $A_n,$
\begin{equation}
\p_n(B) = \frac{\e\la I(R^n\in B) I_{A_n}\ra}{\e\la I_{A_n}\ra}. \
\label{Pndef}
\end{equation}
It is obvious that $\p_n$ is concentrated on the symmetric non-negative definite matrices with off-diagonal elements
now taking values $\{q_1,\ldots, q_{k-1}\}$ and $\p_n$ is invariant under the permutation of replica indices since 
the set $A_n$ is. We will now show that $\p_{n+1}$ restricted to the first $n$ replica coordinates coincides with $\p_n$
and, thus, the sequence $(\p_n)$ defines a law of the infinite overlap array. 
\begin{lemma}\label{Lem1}
For any measurable function $f $ of the overlaps on $n$ replicas, 
\begin{equation}
\e\la f(R^n) I_{A_{n+1}} \ra = (1-p_k) \e\la f(R^n)I_{A_n}\ra.
\label{induction}
\end{equation}
\end{lemma}
\textbf{Proof.}
Notice that $A_n = \{\sigma^1,\ldots,\sigma^n \mbox{ are all different}\}$ by Proposition \ref{Th2} and, therefore,
\begin{equation}
I_{A_{n+1}} = I_{A_n} - \sum_{l\leq n} I_{A_n\cap \{R_{l,n+1} = q_k\}}.
\label{indicator}
\end{equation}
This implies that
$$
\e\la f(R^n) I_{A_{n+1}} \ra = \e\la f(R^n) I_{A_{n}} \ra 
-\sum_{l\leq n} \e\la f(R^n) I_{A_{n}} I(R_{l,n+1} = q_k) \ra.
$$
Using the Ghirlanda-Guerra identities (\ref{GG1}), for each $l\leq n,$
\begin{eqnarray*}
\e\la f(R^n) I_{A_{n}} I(R_{l,n+1} = q_k) \ra
&=&
\frac{p_k}{n} \,\e\la f(R^n) I_{A_{n}} \ra 
+\frac{1}{n}\sum_{l'\not = l}^n \e\la f(R^n) I_{A_{n}} I(R_{l,l'} = q_k) \ra
\\
&=&
\frac{p_k}{n} \,\e\la f(R^n) I_{A_{n}} \ra 
\end{eqnarray*}
since $A_n\subseteq \{R_{l,l'} \not = q_k\}$ and, thus, $I_{A_{n}} I(R_{l,l'} = q_k) = 0.$ 
Adding up over $l\leq n$ finishes the proof.
\qed

First, using (\ref{induction}) inductively for $f\equiv 1$ we get $\e\la I_{A_n}\ra = (1-p_k)^{n-1}$ and then dividing 
(\ref{induction}) by $(1-p_k)^n$ gives
\begin{equation}
 \frac{\e\la f(R^n) I_{A_{n+1}}\ra}{\e\la I_{A_{n+1}}\ra}
 =
 \frac{\e\la f(R^n) I_{A_{n}}\ra}{\e\la I_{A_{n}}\ra}. 
\end{equation}
This means that the family  $(\p_n)$ is consistent and by Kolmogorov's theorem
we can define the distribution of the infinite array with the corresponding marginals given by $\p_n.$
Let us consider an array $Q=(Q_{l,l'})_{l,l'\geq 1}$ with this distribution.

\textbf{Proof of Theorem \ref{ThMain}.}
By construction, $Q$ is a symmetric, non-negative definite and weakly exchangeable array 
with diagonal elements equal to $q_k$ and off-diagonal elements taking values $\{q_1,\ldots,q_{k-1}\}$ with probabilities
\begin{equation}
\p(Q_{1,2} = q_l) = \frac{p_l}{1-p_k}.
\end{equation}
Using the Dovbysh-Sudakov representation for the array $Q$ implies that there exists a random measure $G'$ 
on $H$ such that $Q$ can be generated as 
$$
Q_{l,l'} = \sigma^l\cdot \sigma^{l'}+\delta_{l,l'}(q_k-\sigma^l\cdot\sigma^l)
$$
for an i.i.d. sequence $(\sigma^l)$ from $G'$.
Since $\sigma^l\cdot \sigma^{l'}\in\{q_1,\ldots,q_{k-1}\}$, it is easy to see that the support of $G'$ must be inside the
sphere of radius $\sqrt{q_{k-1}}$ for, otherwise, with positive probability we could sample two points 
$\sigma^1,\sigma^2$ arbitrarily close to a point $\sigma$ such that $\|\sigma\|> \sqrt{q_{k-1}} $ which would contradict
that $\sigma^1\cdot\sigma^2 \leq q_{k-1}$ (see \cite{AA} or \cite{PGG} for details). In particular, the truncated array
$(Q_{l,l'}\wedge q_{k-1})_{l,l'\geq 1}$ can be computed as
\begin{equation}
Q_{l,l'}\wedge q_{k-1} = \sigma^l \cdot\sigma^{l'} + \delta_{l,l'}(q_{k-1}-\sigma^l\cdot\sigma^l)
\label{Qtrunc}
\end{equation}
and it is  non-negative definite as the sum of two non-negative definite arrays. 
If we recall the definition (\ref{Pndef}), the matrix $(Q_{l,l'})_{ l,l'\leq n}$ is obtained by 
sampling $n$ configurations from the measure $G=\sum_{l\geq 1} w_l \delta_{\xi_l}$ conditionally on the
event that these configurations are different. Since with positive probability we can sample
$\xi_{1},\ldots,\xi_{n},$ we must have that the matrix $(\xi_{l}\cdot\xi_{l'} \wedge q_{k-1})_{l,l'\leq n}$ 
is non-negative definite and, therefore, $(\xi_l\cdot\xi_{l'}\wedge q_{k-1})_{l,l'\geq 1}$ is non-negative definite
with probability one. This of course means that $(R_{l,l'}\wedge q_{k-1})_{l,l'\geq 1}$ is also non-negative definite.
Since the function $x\wedge q_{k-1}$ can be approximated by polynomials, the truncated overlap array also satisfies 
the Ghirlanda-Guerra identities and its elements  now take values in $\{q_{1},\ldots,q_{k-1}\}$. One can proceed by induction on $k.$
\qed

\noindent Even though we did not need it in the proof, one can show that the measure $G'$ 
is actually concentrated on the sphere of radius $\sqrt{q_{k-1}}$ by using Proposition \ref{Th2}  and the
following observation.
\begin{lemma}\label{Lem2}
The distribution of $Q$ satisfies the Ghirlanda-Guerra identities,
\begin{equation}
\e f(Q^n) \psi(Q_{1,n+1}) = \frac{1}{n}\, \e f(Q^n)\, \e \psi(Q_{1,2}) + 
\frac{1}{n} \sum_{l=2}^{n} \e f(Q^n) \psi(Q_{1,l}).
\label{GGQ}
\end{equation}
\end{lemma}
\textbf{Proof.} For simplicity of notations let us consider the case of $\psi(x)=x^p.$
Using (\ref{indicator}),
\begin{eqnarray}
\e\la f(R_n) R_{1,n+1}^pI_{A_{n+1}} \ra 
&=& 
\e\la f(R_n) R_{1,n+1}^p I_{A_{n}} \ra 
-\sum_{l\leq n} \e\la f(R_n) R_{1,n+1}^p I_{A_{n}\cap \{R_{l,n+1} = q_k\}} \ra
\nonumber
\\
&=&
\e\la f(R_n) R_{1,n+1}^p I_{A_{n}} \ra 
-\sum_{l\leq n} \e\la f(R_n) R_{1,l}^p I_{A_{n}\cap \{R_{l,n+1} = q_k\}} \ra.
\label{in}
\end{eqnarray}
since $R_{l,n+1}=q_k$ implies that $\sigma^l = \sigma^{n+1}$ and, thus, $R_{1,n+1}= R_{1,l}.$
By the Ghirlanda-Guerra identities, the $l^{th}$ term in the last sum is equal to
\begin{eqnarray*}
&&
\frac{p_k}{n} \,\e\la f(R_n) R_{1,l}^p I_{A_{n}} \ra 
+
\frac{1}{n}\sum_{l'\not = l}^n \e\la f(R_n) R_{1,l}^p I_{A_{n}\cap \{R_{l,l'} = q_k\}} \ra
=
\frac{p_k}{n} \,\e\la f(R_n) R_{1,l}^p I_{A_{n}} \ra 
\end{eqnarray*}
since $A_{n}\cap \{R_{l,l'} = q_k\} = \emptyset.$ 
Similarly,
\begin{eqnarray*}
\e\la f(R_n) R_{1,n+1}^p I_{A_{n}}  \ra
&=&
\frac{1}{n} \,\e\la f(R_n) I_{A_{n}} \ra \e\la R_{1,2}^p\ra
+\frac{1}{n}\sum_{l =2}^n \e\la f(R_n) R_{1,l}^p I_{A_{n} } \ra.
\end{eqnarray*}
Using that $R_{1,1} = q_k$ and combining all the terms in (\ref{in}), $\e\la f(R_n) R_{1,n+1}^p I_{A_{n+1}} \ra$ equals 
\begin{eqnarray*}
&&
\frac{1}{n} \,\e\la f(R_n) I_{A_{n}} \ra 
\bigl( \e\la R_{1,2}^p\ra - q_k^p \,p_k \bigr)
+
\frac{1-p_k}{n}\sum_{l =2}^n \e\la f(R_n) R_{1,l}^p I_{A_{n} } \ra
\\
&&
=\frac{1}{n} \,\e\la f(R_n) I_{A_{n}} \ra  \e\la R_{1,2}^p I_{A_2}\ra 
+
\frac{1}{n}\sum_{l =2}^n \e\la f(R_n) R_{1,l}^p I_{A_{n} } \ra \e\la I_{A_2}\ra.
\end{eqnarray*}
Recalling that $\e\la I_{A_n}\ra = (1-p_k)^{n-1}$ and dividing everything by $(1-p_k)^n$, we get
\begin{eqnarray*}
&&
\frac{\e\la f(R_n) R_{1,n+1}^p I_{A_{n+1}} \ra}{\e\la I_{A_{n+1}}\ra}
=
\frac{1}{n} \frac{\e\la f(R_n) I_{A_{n}} \ra}{\e\la I_{A_n}\ra} \frac{ \e\la R_{1,2}^p I_{A_2}\ra}{\e\la I_{A_2}\ra} 
+
\frac{1}{n}\sum_{l =2}^n \frac{\e\la f(R_n) R_{1,l}^p I_{A_{n} } \ra}{\e\la I_{A_n}\ra}.
\end{eqnarray*}
Comparing with (\ref{Pndef}), this is exactly (\ref{GGQ}).
\qed\\
We would like to point out that the idea of the proof of Theorem \ref{ThMain} suggests the following criterion
of ultrametricity in the general case without the assumption (\ref{finite}). Given $q\in[0,1]$ such that 
$\p(R_{1,2}< q)>0,$ consider the events
\begin{equation}
A_{n,q} = \{R_{l,l'} < q, \, \forall 1\leq l<l'\leq n\}
\end{equation}
and let $\p_{n,q}$ be the distribution of $R^n$ conditionally on $A_{n,q}.$
\smallskip
\begin{theorem}
Under (\ref{GG}), the array $R$ is ultrametric if and only if for any $q$ such that $\p(R_{1,2}< q)>0$ 
and any set $B$ of $3\times 3$ matrices such that $\p_{3,q}(R^3 \in B)>0$ we have $\limsup_{n\to\infty} \p_{n,q}(R^3\in B)>0.$
\end{theorem}
\smallskip
One can check that, in one direction, ultrametricity yields a relationship of the type (\ref{indicator})
which implies the consistency of the sequence $(\p_{n,q})$ as in Lemma \ref{Lem1} and $\p_{n,q}(R^3\in B)=\p_{3,q}(R^3\in B).$
In the other direction, for any $B$ with $\p_{3,q}(R^3\in B)>0$ we can choose the limit $\p_{q}$ over a subsequence of
 $\p_{n,q}$ such that $\p_{q}(R^3\in B)>0$. If ultrametricity fails,  one can make a choice of a subset of non-ultrametric configurations 
 $B$ and $q$ that will lead to contradiction with the Dovbysh-Sudakov representation for $\p_q$.


\begin{thebibliography}{99}

\bibitem{AC} Aizenman, M., Contucci, P. (1998)
On the stability of the quenched state in mean-field spin-glass models. 
{\it J. Statist. Phys.} \textbf{92}, no. 5-6, 765--783. 

\bibitem{AA} Arguin, L.-P., Aizenman, M. (2009)
On the structure of quasi-stationary competing particles systems. 
{\it Ann. Probab.}, \textbf{37}, no. 3, 1080-1113. 

\bibitem{DS} Dovbysh, L. N., Sudakov, V. N. (1982)  Gram-de Finetti matrices. 
{\it Zap. Nauchn. Sem. Leningrad. Otdel. Mat. Inst. Steklov.}  \textbf{119}, 77-86. 

\bibitem{GG} Ghirlanda, S., Guerra, F. (1998) 
General properties of overlap probability distributions in disordered spin systems. 
Towards Parisi ultrametricity.  {\it J. Phys. A}  \textbf{31}, no. 46, 9149-9155.

\bibitem{PGG} Panchenko, D. (2010) A connection between Ghirlanda-Guerra identities and ultrametricity. 
{\it Ann. of Probab.}, \textbf{38}, no. 1, 327-347. 

\bibitem{PDS} Panchenko, D. (2010) On the Dovbysh-Sudakov representation result.
{\it Electron. Comm. in Probab.},  \textbf{15}, 330-338. 

\bibitem{PGGmixed} Panchenko, D. (2010) The Ghirlanda-Guerra identities for mixed $p$-spin model.
{\it C.R.Acad.Sci.Paris, Ser. I}, \textbf{348}, 189-192. 

\bibitem{SG} Talagrand, M. (2003)   Spin Glasses: a Challenge for Mathematicians.  
Ergebnisse der Mathematik und ihrer Grenzgebiete. 3. Folge A Series of Modern Surveys in Mathematics, Vol. 43.
Springer-Verlag.

\bibitem{Tal-New} Talagrand, M. (2010) Construction of pure states in mean-field models
for spin glasses. {\it Probab. Theory Relat. Fields}., \textbf{148}, no. 3-4, 601-643.

\bibitem{SG2} Talagrand, M. (2011) Mean-Field Models for Spin Glasses. 
Ergebnisse der Mathematik und ihrer Grenzgebiete. 3. Folge A Series of Modern Surveys in Mathematics, Vol. 54, 55.
Springer-Verlag.


\end{thebibliography}
\end{document}